\documentclass[%
amsmath,amssymb,
preprint,%
]{revtex4-1}

\usepackage{graphicx}% Include figure files
\usepackage{dcolumn}% Align table columns on decimal point
\usepackage{bm}% bold math
\usepackage[utf8]{inputenc}
\usepackage[T1]{fontenc}
\usepackage{mathptmx}
\usepackage{etoolbox}
\usepackage{amssymb,amsmath,amsthm}
\usepackage{mathrsfs}
\usepackage{hyperref}
\usepackage{bm}
\usepackage{color}
\usepackage{makecell}
\usepackage{float}
\usepackage{overpic}
\usepackage{subfig}
\usepackage{threeparttable}
\newtheorem{theorem}{Theorem}[section]

\newtheorem{lemma}[theorem]{Lemma}

\newtheorem{remark}[theorem]{Remark}

\makeatletter
\def\@email#1#2{%
 \endgroup
 \patchcmd{\titleblock@produce}
  {\frontmatter@RRAPformat}
  {\frontmatter@RRAPformat{\produce@RRAP{*#1\href{mailto:#2}{#2}}}\frontmatter@RRAPformat}
  {}{}
}%
\makeatother
\begin{document}

\preprint{ }

\title{Equivariant Hopf bifurcation arising in circular-distributed predator-prey interaction with taxis}

\author{Yaqi Chen}
\affiliation{Department of Mathematics, Harbin Institute of Technology, Weihai, Shandong 264209, P.R.China.}

\author{Xianyi Zeng}
\affiliation{Department of Mathematics, Lehigh University, Bethlehem, PA 18015, United States.}

\author{Ben Niu~*}%
 \email{Corresponding author: niu@hit.edu.cn}
\affiliation{Department of Mathematics, Harbin Institute of Technology, Weihai, Shandong 264209, P.R.China. %\\This line break forced with \textbackslash\textbackslash
}%

\date{\today}% It is always \today, today,

\begin{abstract}
In this paper, we study the Rosenzweig-MacArthur predator-prey model with predator-taxis and time delay defined on a disk. Theoretically, we studied the equivariant Hopf bifurcation around the positive constant steady-state solution. Standing and rotating waves have been investigated through the theory of isotropic subgroups and Lyapunov-Schmidt reduction. The existence conditions, the formula for the periodic direction and the periodic variation of bifurcation periodic solutions are obtained. Numerically, we select appropriate parameters and conduct numerical simulations to illustrate the theoretical results and reveal quite complicated dynamics on the disk. Different types of rotating and standing waves, as well as more complex spatiotemporal patterns with random initial values, are new dynamic phenomena that do not occur in one-dimensional intervals.
~
\end{abstract}

\maketitle

\section{Introduction}\label{sec1}
In order to survive and reproduce, organisms will respond to the surrounding environment, such as staying away from the environment that hinders their own growth and approaching the environment that is conducive to them. This phenomenon of directed movement is called chemotaxis and research on it can be traced back to \citep{Keller1971J}. In predator-prey models, prey- or predator-taxis can effectively characterize the predator's directional movement toward habitats with high prey density or the prey's behavior of perceiving predation risks to avoid predators, respectively. Many studies have focused on the global existence, boundedness, or stability of solutions in models with prey-taxis \citep{ChenMX2023J,Choi2019J,He2015J,Jin2017J,Giricheva2019J,Wu2016J} or predator-taxis \citep{ChenWu2023J,Wu2018J,Yoon2022J}. Besides, based on the Hopf bifurcation theory for quasilinear reaction-diffusion systems \citep{Amann1991C}, the existence of center manifolds \citep{Simonett1995J}, and the principle of linearized stability \citep{Drangeid1989J}, there are many studies exploring patterns induced by taxis and bifurcations, through center manifold reduction  \citep{Dong2023J,Ma2022J,Shi2022J,Zuo2021J} or Lyapunov-Schmidt reduction method \citep{Gao2022J,Liu2022J,Qiu2020J}.

The above studies mostly focus on one-dimensional intervals. Nevertheless, in real life, the interaction between predators and prey frequently occurs in two-dimensional spatial domains. Examples include lakes \citep{Maciej2006J}, irregularly shaped protected areas \citep{Khan2003J}, and circular Petri dishes, which are often chosen as experimental arenas for many studies on the interaction between miniature predators and prey \citep{Banfield2016J,Meyer2020J}. Establishing a predator-prey model on a circular domain can better depict real-world situations. Hence, we consider the Rosenzweig-MacArthur predator-prey model with predator-taxis and time delay on a disk $\mathbb{D}=\{(r, \theta): 0 \leq r \leq R, 0 \leq \theta \leq 2 \pi\}$,
\begin{equation}\label{system R-M}
\left\{\begin{array}{lll}
u_t=d_1\Delta_{r\theta}u+\chi\nabla_{r\theta}\cdot(u\nabla_{r\theta} v)+u\left(1-\frac{u}{K}\right)-\frac{\alpha uv}{u+1}, & (r,\theta)\in\mathbb{D},~t>0,\\
v_t=d_2\Delta_{r\theta}v-dv+\frac{\alpha u_{\tau}v}{u_{\tau}+1}, & (r,\theta)\in\mathbb{D},~t>0,\\
\partial_r u(\cdot,R,\theta)=\partial_r v(\cdot,R,\theta)=0, & \theta \in [0,2\pi),
\end{array}\right.
\end{equation}
where $u=u(r,\theta,t),~v=v(r,\theta,t)$ represent the density of prey and predator at radius $r$, angle $\theta$ and time $t$, respectively, and $u_{\tau}=u(r,\theta,t-\tau)$. The parameters $d_1,~d_2,~\chi,~K,~\alpha,~d$ are all positive and defined in \citep{Shi2022J}, in which the authors investigated the same model on a one-dimensional interval $[0,l\pi]$.

With the aid of the phase portraits of the normal form with $O(2)$ symmetry \citep{Gils1986J}, the equivariant bifurcation theory \citep{Guo2022J,Qu2023J}, and various wave solutions induced by equivariant bifurcations \citep{Chen2023J,Chen2023J2}, we plan to further explore the joint effect of time delay, predator-taxis and $O(2)$ symmetry.  Compared to results in \citep{Shi2022J}, spatially inhomogeneous periodic solutions on a disk, including standing and rotating waves are found.
The new wave solutions generated in circular domains have greater practical significance and can offer enhanced insights into characterizing population dynamics within high-dimensional spatial domains. In fact, the disk is a very simple abstraction of a two-dimensional compact simple-connected region, where many methods or ideas can be tested on it before applying to more complicated domains. The research methodologies employed to analyze the properties of these wave solutions in this article can also be further extended to more complex domains.

The rest of the paper is organized as follows. In Section \ref{sec2},  we analyze the equivariant Hopf bifurcation around the positive steady state through Lyapunov-Schmidt reduction. In Section \ref{sec4}, we simulated spatially inhomogeneous periodic solutions, including standing waves, rotating waves, and more complex periodic spatiotemporal patterns with random initial values.

\section{Equivariant Hopf bifurcation}\label{sec2}
In this section, we will study the constant steady states of (\ref{system R-M}) through a linearization analysis and analyze the equivariant Hopf bifurcation around the steady state through Lyapunov-Schmidt reduction \citep{Liu2022J,Guo2022J,Qu2023J,Guo2013M}. Firstly, we need to introduce some notations in Table \ref{Notations}.

\begin{table}
\caption{Notations}
\label{Notations}
\setlength{\tabcolsep}{0.02mm}\renewcommand{\arraystretch}{1.25}{
\begin{tabular}{ccccccc}
 \hline
 Symbols                  &  Descriptions\\
    \hline
     $\mathbb{Y}={L}^2(\mathbb{D})$  &  the Lebesgue space of integrable functions defined on $\mathbb{D}$                                                            \\
   $\mathbb{H}_0^1(\mathbb{D})$ & $\mathbb{H}_0^1(\mathbb{D})=\{(u,v)\in \mathbb{H}^1(\mathbb{D})|\partial_r u=\partial_r v=0,~ \theta \in [0,2\pi)\}$\\
   $\mathbb{X}$   &         $\mathbb{X}=\mathbb{H}^2(\mathbb{D}) \cap \mathbb{H}_0^1(\mathbb{D})$ \\
    $\mathbb{X}_{\mathbb{C}},~\mathbb{Y}_{\mathbb{C}}$ & the complexification  of $\mathbb{X}$ and $\mathbb{Y} $\\
   $C([-\tau,0],\mathbb{X}_{\mathbb{C}})$  & the Banach space of continuous mappings form $[-\tau,0]$ to $\mathbb{X}_{\mathbb{C}}$ with  supremum norm \\
   $\lambda_{p},~\lambda_{nm}$ \citep{Chen2023J,Chen2023J2}& the eigenvalues of $-\Delta_{r\theta}$ on $\mathbb{D}$ with homogeneous Neumann boundary conditions\\
   $\hat{\phi}_{p}^c,~\hat{\phi}_{nm}^c,~\hat{\phi}_{nm}^s$ \citep{Chen2023J,Chen2023J2}& the normalised eigenfunctions corresponding to $\lambda_{p},~\lambda_{nm}$\\
   \hline
\end{tabular}
}
\end{table}

In \citep{Shi2022J}, the authors obtained that the system defined on the interval $[0,l\pi]$ undergoes spatially inhomogeneous periodic solution at $\tau=\tau_n^k$,
around the unique positive constant steady state $(u^*,v^*)=\left(\frac{d}{\alpha-d},\frac{(K-u^*)(1+u^*)}{K\alpha}\right)$, under assumptions $\bm{\mathrm{(H_1)}}$-$\bm{\mathrm{(H_3)}}$.
\begin{itemize}
\item [$\bm{\mathrm{(H_1)}}$]
$0<K<1+2u^*,~\alpha>d,~0<u^*<K,~\chi>0.$
\item [$\bm{\mathrm{(H_2)}}$]
$C_{nm}-B_{nm}<0$.
\item [$\bm{\mathrm{(H_3)}}$]
$\chi>\chi_*$.
\end{itemize}
In fact, for $U=(u,v)^{\mathrm{T}},~U_{\tau}=(u_{\tau},v_{\tau})^{\mathrm{T}}$, the linearized system of system \eqref{system R-M} around $(u^*,v^*)$ is
\begin{equation}\label{the linearized system}
\frac{\partial U}{\partial t}=\mathscr{A}_0U+\mathscr{A}_{\tau}U_{\tau},
\end{equation}
where
$$
\mathscr{A}_0= \left(\begin{array}{cc}
d_1\Delta_{r\theta}+a_{11} & \chi u^* \Delta_{r\theta}-d\\
0 & d_2\Delta_{r\theta}
\end{array}\right),~
\mathscr{A}_{\tau}=\left(\begin{array}{cc}
0 & 0\\
a_{21} & 0
\end{array}\right),
$$
with $a_{11}=1-\frac{2u^*}{K}-\frac{\alpha v^*}{(u^*+1)^2},~a_{21}=\frac{\alpha v^*}{(u^*+1)^2}$.
On a disk, the linearized system \eqref{the linearized system} restricted to $\mathrm{span}\left\{\hat{\phi}_{p}^c,~\hat{\phi}_{nm}^c,~\hat{\phi}_{nm}^s\right\}$ is equivalent to a sequence of  functional differential equations (FDEs) on $\mathbb{R}$, with a sequence of characteristic equations. Due to the symmetry, some of them have the following form
\begin{equation}\label{Gamma1}
\Gamma(\gamma,\lambda_{p},\chi)=\gamma^2+A_{p}\gamma+B_{p}\mathrm{e}^{-\gamma \tau}+C_{p}=0,~p=0,1,2,\cdots,
\end{equation}
and the others are in the following form \citep{Chen2023J,Chen2023J2}
\begin{equation}\label{Gamma2}
\Gamma(\gamma,\lambda_{nm},\chi)=\left[\gamma^2+A_{nm}\gamma+B_{nm}\mathrm{e}^{-\gamma \tau}+C_{nm}\right]^2=0,~n>0,~m>0,
\end{equation}
with $A_{nm}=(d_1+d_2)\lambda_{nm}-a_{11},~B_{nm}=a_{21}(\chi u^*\lambda_{nm}+d),~C_{nm}=(d_1\lambda_{nm}-a_{11})d_2\lambda_{nm}$, and the expressions for $A_{p},~B_{p},~C_{p}$ can be obtained analogously. In this article, we are more concerned about the phenomenon caused by symmetry, so we will study the problem in the case of $n>0,~m>0$.
To distinguish the symbols with \citep{Shi2022J}, we record the critical value of the bifurcation parameter $\tau$ as $\tau_{{nm}}^k$, and
 Eq.(\ref{Gamma2}) has a pair of repeated purely imaginary roots $\gamma =\pm \mathrm{i}\omega_{*}(\omega_{*}>0)$, with $\omega_*=\sqrt{\left({(2C_{nm}-A_{nm})+\sqrt{(2C_{nm}-A_{nm})^2-4(C_{nm}-B_{nm})}}\right)/{2}}$, under the same assumptions $\bm{\mathrm{(H_1)}}$-$\bm{\mathrm{(H_3)}}$. Besides, the
transversality conditions hold, i.e. $\mathrm{Re}\{\gamma^{\prime}({\tau_{{nm}}^k})\} \ne 0$. That is to say,  system \eqref{system R-M} will generate equivariant Hopf bifurcation at $\tau=\tau_{{nm}}^k$ on the disk,  by \citep{Shi2022J,Chen2023J,Guo2013M}.

In the following content, we will study the equivariant Hopf bifurcation around $(u^*,v^*)$ and find a periodic solution with a period near $T=\frac{2\pi}{\omega_*}$. For preparation work,
%Note that
%\begin{equation}\label{M}
%M(\gamma,\tau,\chi)=\left[\begin{array}{cc}
%-d_1\lambda_{nm}+a_{11}-\gamma & -\chi u^*\lambda_{nm}-d\\
%a_{21}\mathrm{e}^{-\gamma \tau} & -d_2\lambda_{nm}-\gamma
%\end{array}\right]^2.
%\end{equation}
%Obviously, we have $M(\mathrm{i}\omega_{*},\tau_{\lambda_{nm}}^k,\chi)V_1=0$ and $M^{\mathrm{T}}(-\mathrm{i}\omega_{*},\tau_{\lambda_{nm}}^k,\chi)V_2=0$, where
%
%Obviously, $V_2\cdot V_1=0$ and $V_2\cdot \bar{V}_1=1$.
%From \citep{Shi2022J}, for all $\tau$ near $\tau_{\lambda_{nm}}^k$, we have $M(\gamma,\tau, \chi)=0$. Thus, there exists a continuously differentiable mapping $V(\tau)$ such that $V(\tau_{\lambda_{nm}}^k)=V_1$ and $M(\gamma,\tau, \chi)V(\tau)=0$. Differentiating $M(\gamma,\tau, \chi)V(\tau)=0$ with respect to $\tau$ and letting $\tau=\tau_{\lambda_{nm}}^k$, we obtain
%$$
%\left[M_{\tau}(\mathrm{i}\omega_{*},\tau_{\lambda_{nm}}^k,\chi)+{\gamma}^{\prime}(\tau_{\lambda_{nm}}^k)M_{\gamma}(\mathrm{i}\omega_{*},\tau_{\lambda_{nm}}^k,\chi)\right]V_1
%+M(\mathrm{i}\omega_{*},\tau_{\lambda_{nm}}^k,\chi)V^{\prime}(\tau)=0.
%$$
%Note that $M^{\mathrm{T}}(-\mathrm{i}\omega_{*},\tau_{\lambda_{nm}}^k,\chi)V_2=0$, thus, $\bar{V}_2\cdot M_{\tau}(\mathrm{i}\omega_{*},\tau_{\lambda_{nm}}^k,\chi)V_1={\gamma}^{\prime}(\tau_{\lambda_{nm}}^k)$.
let
$$
\begin{aligned}
\mathbb{C}_T~(\mathbb{C}_T^1)=&\left\{{F}(t,\cdot,\cdot): \mathbb{R}\rightarrow \mathbb{X}_{\mathbb{C}}|{F}(t+T,\cdot,\cdot)={F}(t,\cdot,\cdot)\right.\\
&\left.{\rm{~and~is ~continuous~(differentiable)~with~respect~to}}~t\right\}.
\end{aligned}
$$
For ${F}(t)\in \mathbb{C}_T~(\mathbb{C}_T^1)$, note that
$
\|{F}(t)\|_0=\max_{t\in[0,T]}\{\|{F}(t)\|_{\mathbb{X}_{\mathbb{C}}^2}\}$ $\left(\|{F}(t)\|_1=\max\{\|{F}(t)\|_0,\|{F}^{\prime}(t)\|_0\}\right)$.
Then, $\mathbb{C}_T$ and $\mathbb{C}_T^1$ are Banach spaces with norms $\|{F}(t)\|_0$ and $\|{F}(t)\|_1$, respectively. In fact, $\mathbb{C}_T$ is a Banach representation of the group $\mathbb{S}^1$, i.e.
$$
\kappa \cdot{F}(t)={F}(t+\kappa),~\kappa \in \mathbb{S}^1, {F} \in\mathbb{C}_{2\pi}.
$$
For all ${F},~{G}$ in $\mathbb{C}_T$, define the inner product $(\cdot,\cdot)_T:\mathbb{C}_T\times\mathbb{C}_T\rightarrow\mathbb{R}$ as
$
  ({F},{G})_T=\frac{1}{T}\int_0^T\langle{F}(t),{G}(t)\rangle\mathrm{d}t,
$
where $\langle u(r,\theta),v(r,\theta)\rangle=\iint_{\mathbb{D}}r u(r,\theta) \bar{v}(r,\theta) \mathrm{d} r\mathrm{d}\theta$ is the inner product weighted $r$ on $\mathbb{Y}_{\mathbb{C}}$.

Normalizing of the period, let $\rho\in(-1,1)$ and $U(t)=\left(u((1+\rho)t),v((1+\rho)t)\right)^{\mathrm{T}}$, then system (\ref{system R-M}) can be transformed into
\begin{equation}\label{system R-M U}
  (1+\rho)\frac{\mathrm{d}U(t)}{\mathrm{d}t}=\mathscr{A}_0U(t)+\mathscr{A}_{\tau}U(t-(1+\rho)\tau)+\mathscr{H}(U_t,\tau),
\end{equation}
where
$
\mathscr{H}(\phi,\tau)=\frac{1}{2}\mathscr{B}(\phi,\phi)+\frac{1}{6}\mathscr{C}(\phi,\phi,\phi)+o(\|\phi\|^3)$, for $\phi \in \mathbb{X}_{\mathbb{C}}^2.
$
The above transformation simplifies the study of periodic solutions with time $T$ and its similar times of system \eqref{system R-M}. Finding $\frac{T}{1+\rho}$-periodic solutions of system (\ref{system R-M}) is equivalent to solving $\mathscr{F}(U,\tau,\rho)=0$, where $\mathscr{F}:\mathbb{C}_{T}^1\times\mathbb{R}^2\rightarrow\mathbb{C}_{T}$ is defined by
\begin{equation}\label{huaF}
\mathscr{F}(U,\tau,\rho)=-(1+\rho)\frac{\mathrm{d}U(t)}{\mathrm{d}t}+\mathscr{A}_0U(t)+\mathscr{A}_{\tau}U(t-(1+\rho)\tau)+\mathscr{H}(U_t,\tau).
\end{equation}
Obviously, $\mathscr{F}$ is $O(2)\times \mathbb{S}^1$-equivariant, i.e. for $(\delta,\kappa)\in O(2)\times \mathbb{S}^1$,
$$
(\delta,\kappa)\mathscr{F}(U,\tau,\rho)=\mathscr{F}\left( U\left(r,\theta+\frac{\delta}{n},t+\frac{\kappa}{\omega_*(1+\rho)}\right),\tau,\rho\right),~ (\delta~stands~for~a~rotation)
$$
or
$$
(\delta,\kappa)\mathscr{F}(U,\tau,\rho)=\mathscr{F}\left( U\left(r,-\theta,t+\frac{\kappa}{\omega_*(1+\rho)}\right),\tau,\rho\right), ~(\delta~stands~for~a~reflection)
$$

Let $\mathscr{L}_{\tau_{nm}^k}: \mathbb{X}_{\mathbb{C}} \rightarrow \mathbb{Y}_{\mathbb{C}}$ be the first derivative of $\mathscr{F}$ with respect to $U(t)$ at $(U,\tau,\rho)=(0,\tau_{{nm}}^k,0)$. Clearly, the elements of ${\rm{Ker}}\mathscr{L}_{\tau_{{nm}}^k}$ correspond to the $T$-periodic solutions of the linearized equation $\mathscr{L}_{\tau_{{nm}}^k}U=0$. Let $\mathscr{L}_{\tau_{{nm}}^k}^*$ be the adjoint operator of $\mathscr{L}_{\tau_{{nm}}^k}$ with $(\cdot,\cdot)$ that is a kind of bilinear form, defined in \citep{Chen2023J,Guo2013M}. Then we have the following results.
\begin{lemma}
Let $$V_1=(1,
\frac{a_{21}\mathrm{e}^{\mathrm{i}\omega_*\tau_{{nm}}^k}}{d_2\lambda_{nm}+\mathrm{i}\omega_*})^{\mathrm{T}},
~V_2=\frac{-d_2\lambda_{nm}+\mathrm{i}\omega_*}{(d_2\lambda_{nm}+\mathrm{i}\omega_*)\mathrm{e}^{\mathrm{i}\omega_*\tau_{{nm}}^k}-d_2\lambda_{nm}+\mathrm{i}\omega_*}(
1,
\frac{d_2\lambda_{nm}+\mathrm{i}\omega_*}{-a_{21}\mathrm{e}^{\mathrm{i}\omega_*\tau_{{nm}}^k}}).$$
For $\vartheta \in [-\tau,0]$, define
$$\{ \varphi_1(\vartheta),\varphi_2(\vartheta),\varphi_3(\vartheta),\varphi_4(\vartheta)\}=\{\mathrm{e}^{\mathrm{i}\omega_* \vartheta} V_1 \hat{\phi}_{nm}^c,\mathrm{e}^{-\mathrm{i}\omega_*  \vartheta} \bar{V}_1 \hat{\phi}_{nm}^c,\mathrm{e}^{\mathrm{i}\omega_* \vartheta} V_1 \hat{\phi}_{nm}^s,\mathrm{e}^{-\mathrm{i}\omega_*  \vartheta} \bar{V}_1 \hat{\phi}_{nm}^s\},$$
and
$$\{\varphi_1^*(\vartheta),\varphi_2^*(\vartheta),\varphi_3^*(\vartheta),\varphi_4^*(\vartheta)\}=\{\mathrm{e}^{\mathrm{i}\omega_* \vartheta} V_2 \hat{\phi}_{nm}^c,\mathrm{e}^{-\mathrm{i}\omega_*  \vartheta} \bar{V}_2 \hat{\phi}_{nm}^c,\mathrm{e}^{\mathrm{i}\omega_* \vartheta} V_2 \hat{\phi}_{nm}^s,\mathrm{e}^{-\mathrm{i}\omega_*  \vartheta} \bar{V}_2 \hat{\phi}_{nm}^s\}.$$
Then \\
$\mathrm{(i)}$~${\rm{Ker}}\mathscr{L}_{\tau_{{nm}}^k}$ is spanned by $\left\{\varphi_1,\varphi_2,\varphi_3,\varphi_4\right\}$ and ${\rm{Ker}}\mathscr{L}_{\tau_{{nm}}^k}^*$ is spanned by $\left\{\varphi_1^*,\varphi_2^*,\varphi_3^*,\varphi_4^*\right\}$.\\
$\mathrm{(ii)}$~$(\varphi_j^*,\varphi_j)=1$, for $j=1,2,3,4$, and $(\varphi_i^*,\varphi_j)=0$, for $i \ne j$.
\end{lemma}

Obviously, $\mathscr{L}_{\tau_{{nm}}^k}$ is a Fredholm operator with index 0 \citep{Gao2022J,Qu2020J}, and according to \citep{Guo2013M}, we can get a space decomposition by $O(2)\times\mathbb{S}^1$-invatiant subspaces of $\mathbb{C}_T$ :
$$
\mathbb{C}_T={\rm{Ker}}\mathscr{L}_{\tau_{{nm}}^k}\oplus{\rm{Ran}}\mathscr{L}_{\tau_{{nm}}^k},~
\mathbb{C}_T^1={\rm{Ker}}\mathscr{L}_{\tau_{{nm}}^k}\oplus\mathbb{W},
$$
where $\mathbb{W}=({\rm{Ker}}\mathscr{L}_{\tau_{{nm}}^k}^*)^{\perp}\cap \mathbb{C}_T^1$. Next, we aim to obtain the bifurcation mapping corresponding to $\mathscr{F}(U,\tau,\rho)=0$ through Lyapunov-Schmidt reduction. Define two project operators $P:\mathbb{C}_T \rightarrow {\rm{Ran}}\mathscr{L}_{\tau_{{nm}}^k}$ and $I-P:\mathbb{C}_T \rightarrow {\rm{Ker}}\mathscr{L}_{\tau_{{nm}}^k}$. Then, $P$ and $I-P$ are $O(2)\times\mathbb{S}^1$-equivariant. According to the above space decomposition, $\mathscr{F}(U,\tau,\rho)=0$ can be decomposed into the following two equivalent equations:
$$
\begin{aligned}
&\mathscr{F}_1(V,W,\tau,\rho)=P\mathscr{F}(V+W,\tau,\rho)=0,\\
&\mathscr{F}_2(V,W,\tau,\rho)=(I-P)\mathscr{F}(V+W,\tau,\rho)=0,
\end{aligned}
$$
where $U=V+W,~V\in{\rm{Ker}}\mathscr{L}_{\tau_{{nm}}^k},~W\in\mathbb{W}$. Due to $\mathscr{F}_1(0,0,\tau_{{nm}}^k,0)=0,~D_{W}\mathscr{F}_2(0,0,\tau,0)=P\mathscr{L}_{\tau_{{nm}}^k}=\mathscr{L}_{\tau_{{nm}}^k}$, and $\left.\mathscr{L}_{\tau_{{nm}}^k}\right|_{\mathbb{W}}$ is invertible, by the implicit function theorem, we get a continuously differentiable $O(2)\times \mathbb{S}^1$-equivariant map $\mathscr{W}:\mathscr{L}_{\tau_{{nm}}^k}\times\mathbb{R}^2\rightarrow \mathbb{W}$ and satisfies $\mathscr{W}(0,\tau_{{nm}}^k,0)=0$ and $\mathscr{F}_1(V+\mathscr{W}(V,\tau,\rho),\tau,\rho)=0$. Replacing $W$ in $\mathscr{F}_2(V,W,\tau,\rho)$ by $\mathscr{W}$, we have
\begin{equation}\label{huaF2W}
\mathscr{F}_2^{\mathscr{W}}(V,\tau,\rho)=(I-P)\mathscr{F}(V+\mathscr{W}(V,\tau,\rho),\tau,\rho)=0.
\end{equation}
It is easy to see $\mathscr{F}_2^{\mathscr{W}}$ is $O(2)\times \mathbb{S}^1$-equivariant and satisfies $\mathscr{F}_2^{\mathscr{W}}(0,\tau_{{nm}}^k,0)=0$ and $D_V\mathscr{F}_2^{\mathscr{W}}(0,\tau_{{nm}}^k,0)=0$.
Therefore, by \citep{Guo2013M}, the periodic solution with time period close to $\frac{2\pi}{\omega_*}$ of system \eqref{system R-M} corresponds to the solution of \eqref{huaF2W}.

%For each $V\in {\rm{Ker}}\mathfrak{L}_{\tau_{\lambda_{nm}}^k}$, we have $V=z_1\varphi_1+z_2\varphi_2+z_3\varphi_3+z_4\varphi_4$, where $z_j=(\varphi_j^*,V),~j=1,2,3,4$. Let
%\begin{equation}\label{B}
%\left\{\begin{array}{cc}
%\mathfrak{B}^1(\varepsilon,\tau,\rho)=(\varphi_1^*,\mathfrak{F}_2^{\mathscr{W}}(z_1\varphi_1+z_2\varphi_2+z_3\varphi_3+z_4\varphi_4,\tau,\rho))\\
%\mathfrak{B}^2(\varepsilon,\tau,\rho)=(\varphi_3^*,\mathfrak{F}_2^{\mathscr{W}}(z_1\varphi_1+z_2\varphi_2+z_3\varphi_3+z_4\varphi_4,\tau,\rho))
%\end{array}\right.
%\end{equation}
%It follows that $\mathfrak{B}^1$ and  $\mathfrak{B}^2$ are $\mathbb{S}^1$-equivariant and satisfies $\mathfrak{B}^1_{z_j}(0,\tau_{\lambda_{nm}}^k,0)=0$ and $\mathfrak{B}^2_{z_j}(0,\tau_{\lambda_{nm}}^k,0)=0,~j=1,2,3,4$. By \citep{Golubitsky1989M}, there exist two functions $\mathfrak{g}^1$ and

However, finding periodic solutions of system \eqref{system R-M} is closely related to the $O(2)\times \mathbb{S}^1$ symmetry, which often requires a subspace with simple eigenvalues to study.
Therefore, we can discuss based on the isotropy subgroup theory. %$O(2)\times\mathbb{S}^1$ is isomorphic to $\mathbb{C}^2$, and the action of $O(2)\times\mathbb{S}^1$ on $\mathbb{C}^2$ is given by
%$$
%\delta_{\bar{\theta}}\cdot(u,v)=(u,v)M_{n\bar{\theta}},~\delta_{z}\cdot(u,v)=(u,v)M_{z},~\kappa\cdot(u,v)=(u,v)\mathrm{e}^{\mathrm{i}\kappa},
%$$
%where $(u,v)\in \mathbb{C}^2,~\bar{\theta} \in \mathbb{R},~\kappa \in \mathbb{S}^1$, $\delta_{\bar{\theta}}$ represents rotation through angle $\bar{\theta}$ and $\delta_{z}$ represents the reflection, with
%$$
%M_{\bar{\theta}}=\left[\begin{array}{cc}
%\cos{\bar{\theta}} & -\sin{\bar{\theta}}\\
%\sin{\bar{\theta}} & \cos{\bar{\theta}}
%\end{array}\right],
%~M_{z}=\left[\begin{array}{cc}
%1 & 0\\
%0 & -1
%\end{array}\right].
%$$
By \citep{Qu2023J,Guo2013M,Golubitsky1989M}, the maximal isotropy subgroup of $O(2)\times \mathbb{S}^1$ is $\widetilde{SO}(2)$ and $\mathbb{Z}_2\oplus\mathbb{Z}_2^{\mathbb{C}}$.
%where $\widetilde{SO}(2)$ is generated by $(\delta_{\bar{\theta}/n},\bar{\theta})$, $\mathbb{Z}_2$ is generated by $(\delta_{z},0),~(\delta_{\pi/n}\delta_{z},\pi)$ and $\mathbb{Z}_2^{\mathbb{C}}$ is generated by $(\delta_{\pi/n},\pi)$.
Obviously, the $\widetilde{SO}(2)$-symmetric solutions correspond to rotating wave solutions, and satisfy that
\begin{equation} \label{clockwise-rotating}  (\bar{\theta},\bar{\theta})u(r,\theta,t)=u\left(r,\theta+\frac{\bar{\theta}}{n},t+\frac{\bar{\theta}}{\omega_*}\right),~\bar{\theta} \in \mathbb{S}^1,
\end{equation}
\begin{equation} \label{counterclockwise-rotating}
(\bar{\theta},-\bar{\theta})u(r,\theta,t)=u\left(r,\theta+\frac{\bar{\theta}}{n},t-\frac{\bar{\theta}}{\omega_*}\right),~\bar{\theta} \in \mathbb{S}^1,
\end{equation}
which correspond to the two possible senses of rotation: clockwise and counterclockwise. Note that  $\mathrm{Fix}(\widetilde{SO}(2), \\
{\rm{Ker}}\mathscr{L}_{\tau_{{nm}}^k})=\mathrm{span}\{\psi,\bar{\psi}\}$ with $\psi=\varphi_3,~\bar{\psi}=\varphi_2$ or $\psi=\varphi_1,~\bar{\psi}=\varphi_4$, and $\psi^*=\varphi_3^*$ or $\psi^*=\varphi_1^*$. It is easy to see that $\mathrm{dim}\mathrm{Fix}(\widetilde{SO}(2),{\rm{Ker}}\mathscr{L}_{\tau_{{nm}}^k})=2$.
Analogously, the $\mathbb{Z}_2\oplus\mathbb{Z}_2^{\mathbb{C}}$-symmetric solutions correspond to a standing wave solutions, and satisfy that
\begin{equation}\label{standing form}
\begin{aligned}
 \kappa u(r,\theta,t)=u(r,-\theta,t),&~(\pi,\pi)u(r,\theta,t)=u\left(r,\theta+\frac{\pi}{n},t+\frac{T}{2}\right),\\
\kappa(\pi,\pi)u(r,\theta,t)&=u\left(r,-\theta-\frac{\pi}{n},t+\frac{T}{2}\right).
\end{aligned}
\end{equation}
Note that $\mathrm{Fix}(\mathbb{Z}_2\oplus\mathbb{Z}_2^{\mathbb{C}},{\rm{Ker}}\mathscr{L}_{\tau_{{nm}}^k})=\mathrm{span}\{\psi,\bar{\psi}\}$ with $\psi=\varphi_1+\varphi_3,~\bar{\psi}=\varphi_4+\varphi_2$, and $\psi^*=\varphi_1^*+\varphi_3^*$.

Inspired by \citep{Qu2023J,Guo2013M}, we simplify the problem once again and consider a new restriction mapping $\mathscr{M}:\mathrm{Fix}(\Sigma,{\rm{Ker}}\mathscr{L}_{\tau_{{nm}}^k}) \times \mathbb{R}^2\rightarrow {\rm{Ker}}\mathscr{L}_{\tau_{{nm}}^k},~\Sigma=\widetilde{SO}(2)$ or $\mathbb{Z}_2\oplus\mathbb{Z}_2^{\mathbb{C}}$,
$$
\mathscr{M}(V,\tau,\rho)=\mathscr{F}_2^{\mathscr{W}}(V,\tau,\rho).
$$
Then, $\mathscr{M}(V,\tau,\rho)=0$ is equivalent to the one-dimensional bifurcation map $\mathscr{G}(z,\tau,\rho):\mathbb{C}\times\mathbb{R}^2\rightarrow\mathbb{C}$
\begin{equation}\label{GG}
\mathscr{G}(z,\tau,\rho)=(\psi^*,\mathscr{M}(z\psi+\bar{z}\bar{\psi},\tau,\rho))_T,
\end{equation}
where $z=\langle\psi^*,V\rangle$. Now, our problem is reduced to finding solutions of \eqref{GG}. Clearly, $\mathscr{G}(z,\tau,\rho)$ is $\mathbb{S}^1$-equivariant, and satisfies $\mathscr{G}_z(0,{\tau_{{nm}}^k},0)=0,~\mathscr{G}_{\bar{z}}(0,{\tau_{{nm}}^k},0)=0$. Therefore, there exist two functions ${g}^1,{g}^2:\mathbb{R}^3\rightarrow\mathbb{R}$, such that
\begin{equation}\label{G}
\mathscr{G}(z,\tau,\rho)={g}^1(|z|^2,\tau,\rho)z+{g}^2(|z|^2,\tau,\rho)\mathrm{i}z.
\end{equation}
Since $\mathscr{G}_z(0,{\tau_{{nm}}^k},0)=0$, then ${g}^1(0,{\tau_{{nm}}^k},0)={g}^2(0,{\tau_{{nm}}^k},0)=0$. Let $z=\tilde{\rho}\mathrm{e}^{\mathrm{i}\tilde{\theta}}$, then $\mathscr{G}(z,\tau,\rho)=0$ is equivalent to either solving $\tilde{\rho}=0$ or ${g}^1(\tilde{\rho}^2,\tau,\rho)z=0$ and ${g}^2(\tilde{\rho}^2,\tau,\rho)z=0$. Notice that
$$
\begin{aligned}
&\mathscr{G}_{\tau}(z,{\tau_{{nm}}^k},0)=z\gamma^{\prime}({\tau_{{nm}}^k})+O(|z|^2),\\
&\mathscr{G}_{\rho}(z,{\tau_{{nm}}^k},0)=-\mathrm{i}\omega_*z+O(|z|^2).
\end{aligned}
$$
Therefore,
$$
\mathrm{det}\left[\begin{array}{cc}
{g}^1_{\tau}(0,{\tau_{{nm}}^k},0)  &{g}^1_{\rho}(0,{\tau_{{nm}}^k},0)\\
{g}^2_{\tau}(0,{\tau_{{nm}}^k},0)  &{g}^2_{\rho}(0,{\tau_{{nm}}^k},0)
\end{array}\right]=
\mathrm{det}\left[\begin{array}{cc}
\mathrm{Re}\{\gamma^{\prime}({\tau_{{nm}}^k})\}  &\mathrm{Im}\{\gamma^{\prime}({\tau_{{nm}}^k})\}\\
0  &-\omega_*
\end{array}\right] \ne 0.
$$
The implicit function theorem implies that there exists two functions ${\tau}={\tau}(\tilde{\rho}^2)$ and ${\rho}={\rho}(\tilde{\rho}^2)$ satisfying ${\tau}(0)={\tau_{{nm}}^k}$, ${\rho}(0)=0$ and
\begin{equation}\label{g1g2}
{g}^1(\tilde{\rho}^2,{\tau}(\tilde{\rho}^2),{\rho}(\tilde{\rho}^2))=0,~{g}^2(\tilde{\rho}^2,{\tau}(\tilde{\rho}^2),{\rho}(\tilde{\rho}^2))=0,
\end{equation}
for all sufficient small $\tilde{\rho}$. Therefore, we have the following equivariant Hopf bifurcation theorem for system (\ref{system R-M}).
\begin{theorem}
For every maximal isotropy subgroup $\Sigma \subset O(2)\times \mathbb{S}^1$ in which $\Sigma=\widetilde{SO}(2)$ or $\mathbb{Z}_2\oplus\mathbb{Z}_2^{\mathbb{C}}$,  $\mathrm{Fix}(\Sigma,\mathbb{X})=\{U\in\mathbb{X} | \delta \cdot U=U,~\delta \in \Sigma\}$ is of dimension two, system (\ref{system R-M}) has a branch of periodic solutions emanating from $(u^*,v^*)$ at $\tau={\tau_{{nm}}^k}$, whose spatiotemporal symmetry can be completely characterized by $\Sigma$, corresponding to a branch of clockwise (counterclockwise) rotating waves of the form (\ref{clockwise-rotating}) (form (\ref{counterclockwise-rotating})) or a branch of standing waves of the form (\ref{standing form}).
\end{theorem}

In what follows, we consider the bifurcation direction and the monotonicity of the period. We use
\begin{equation}\label{g21}
 g_{21}=\left(\psi^*,\mathscr{C}(\psi,\psi,\bar{\psi})\right)_T+2\left(\psi^*,\mathscr{B}(\psi,W_{11})\right)_T+\left(\psi^*,\mathscr{B}(\bar{\psi},W_{20})\right)_T,
\end{equation}
to denotes the coefficient of the term $z^2\bar{z}$ in the Taylor expansion of $\mathscr{G}(z,{\tau_{{nm}}^k},0)$. We still need to compute $W_{11}$ and $W_{20}$, which denote the coefficient of the terms $z\bar{z}$ and $z^2$ in the Taylor expansion of $\mathscr{W}(z\psi+\bar{z}\bar{\psi},{\tau_{{nm}}^k},0)$, respectively. In fact,
$$
\begin{aligned}
W_{11}=-\mathscr{L}_{\tau_{{nm}}^k}^{-1}P\mathscr{B}(\psi,\bar{\psi}),~W_{20}=-\mathscr{L}_{\tau_{{nm}}^k}^{-1}P\mathscr{B}(\psi,{\psi}).
\end{aligned}
$$
Namely, $\mathscr{B}(\psi,\bar{\psi}),~\mathscr{B}(\psi,{\psi})\in \mathrm{Ran}\mathscr{L}_{\tau_{{nm}}^k}$, thus, the projection operator $P$ on them acts as the identity. Following the definition of $\mathscr{L}_{\tau_{{nm}}^k}$, we have
$$
\begin{aligned}
&W_{11}=-\left[\begin{array}{cc}
d_1\Delta_{r\theta}+a_{11} & \chi u^* \Delta_{r\theta}-d\\
0 & d_2\Delta_{r\theta}
\end{array}\right]^{-1}\mathscr{B}(\psi,\bar{\psi}),\\
&W_{20}=-\left[\begin{array}{cc}
d_1\Delta_{r\theta}+a_{11} & \chi u^* \Delta_{r\theta}-d\\
a_{21}\mathrm{e}^{-2\mathrm{i}\omega_*\tau_{{nm}}^k} & d_2\Delta_{r\theta}
\end{array}\right]^{-1}\mathscr{B}(\psi,\bar{\psi}).
\end{aligned}
$$
By \citep{Guo2022J,Guo2013M,Faria2000J,Wu1996M}, we summarize the following results.
\begin{theorem}\label{xingzhiLS}
There exists a branch of $\Sigma$-symmetric periodic solutions ($\Sigma=\widetilde{SO}(2)$ or $\mathbb{Z}_2\oplus\mathbb{Z}_2^{\mathbb{C}}$), parameterized by $\tau$, bifurcating from the positive constant steady state of \eqref{system R-M}. Moreover,\\
(i) ${\tau}^{\prime}(0)=\frac{\mathrm{Re}\{g_{21}\}}{\mathrm{Re}\{\gamma^{\prime}(\tau_{{nm}}^k)\}}$ determines the direction of the bifurcation. The bifurcation is supercritical (subcritical), if ${\tau}^{\prime}(0)<0~(>0)$.\\
(ii)  ${\rho}^{\prime}(0)=\frac{\mathrm{Im}\{\gamma^{\prime}(\tau_{{nm}}^k)\bar{g}_{21}\}}{\mathrm{Re}\{\gamma^{\prime}(\tau_{{nm}}^k)\}}$ determines the period of the bifurcation periodic solutions. The period is greater (smaller) than $\frac{2\pi}{\omega_*}$, if ${\rho}^{\prime}(0)>0~(<0)$.
\end{theorem}

\section{Numerical simulations}\label{sec4}

\begin{figure}[h]
\centering
\includegraphics[width=0.55\textwidth]{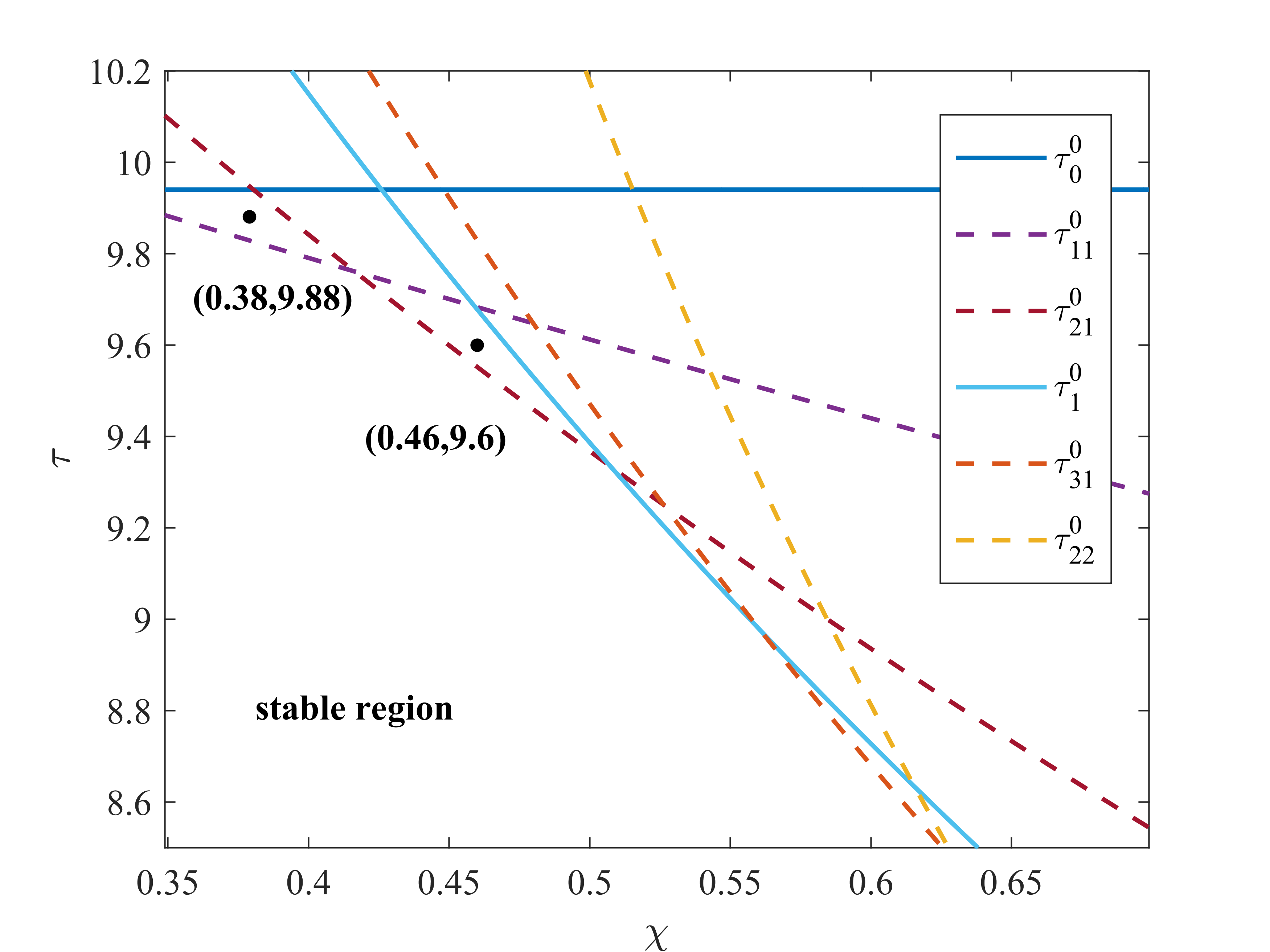}
\caption{Partial bifurcation curves on the $\chi-\tau$ plane. Parameters are $d_1=0.1,~d_2=0.2,~\alpha=1,~K=6,~d=0.1,~R=10$.}
\label{fenzhitutaxis}
\end{figure}

\begin{figure}[t]
\centering
(a)\includegraphics[width=0.47\textwidth]{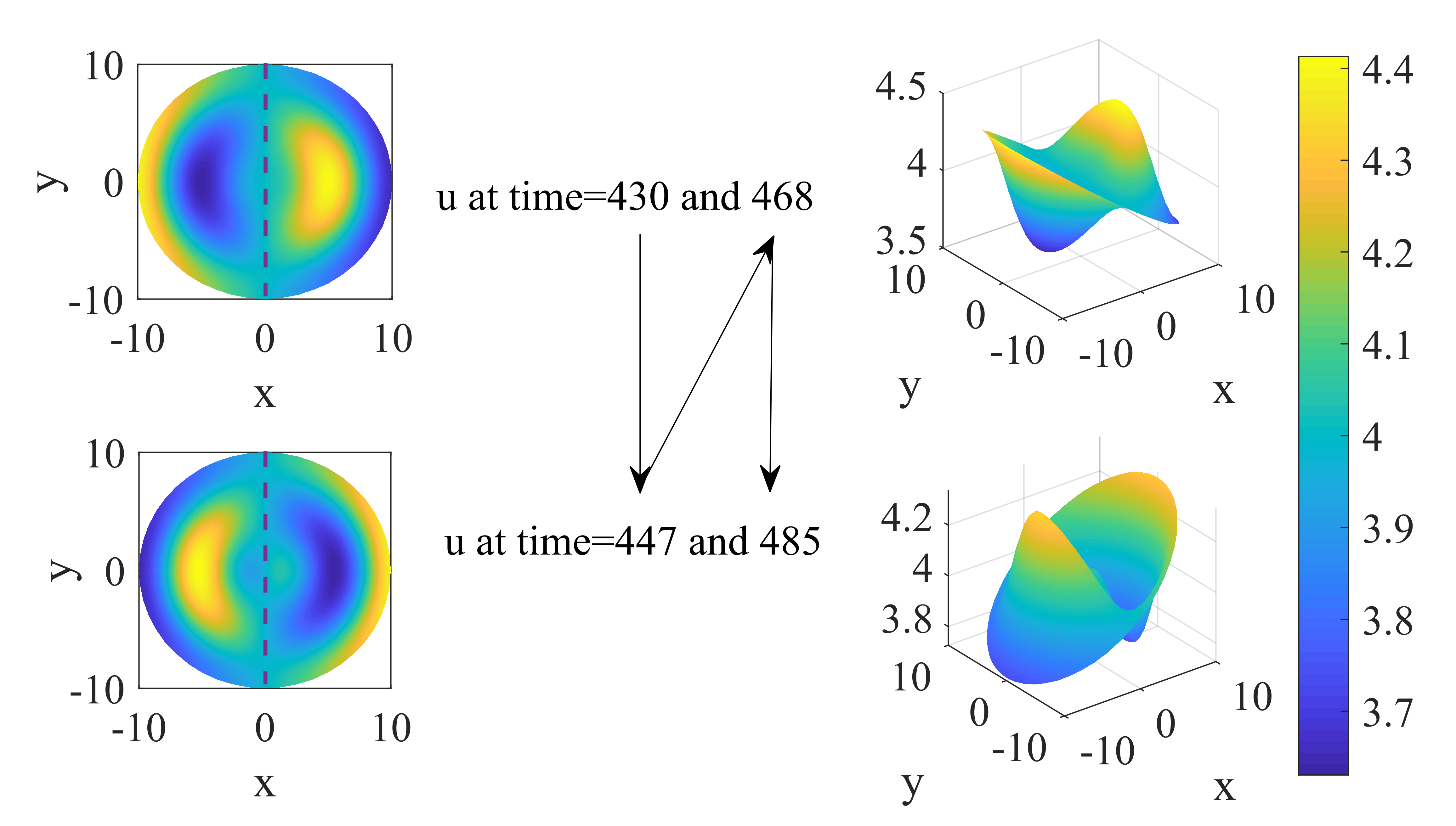}
(b)\includegraphics[width=0.47\textwidth]{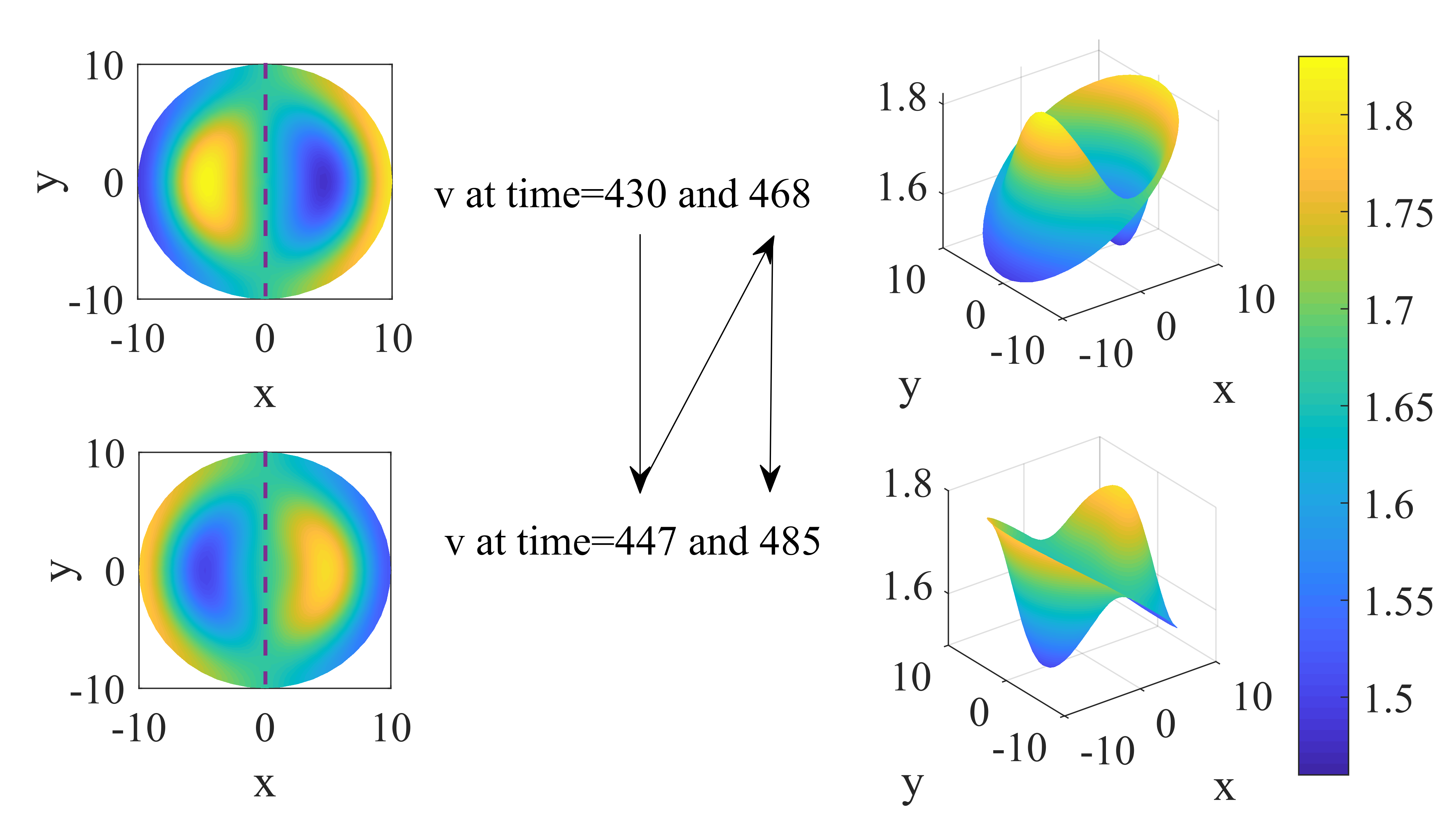}
\caption{System (\ref{system R-M}) generates standing wave that has a fixed axis with $(\chi,\tau)=(0.38,9.88)$.
Initial values are $u(t,r,\theta)=u^*(1+0.1\cdot\cos{t}\cdot\cos{\frac{2\pi r}{R}}\cdot\cos{\theta}),~v(t,r,\theta)=v^*(1+0.1\cdot\cos{t}\cdot\cos{\frac{2\pi r}{R}}\cdot\cos{\theta}),~t\in[-\tau,0)$. (a): u; (b): v.}
\label{S-2}
\end{figure}

\begin{figure}[t]
\centering
(a)\includegraphics[width=0.47\textwidth]{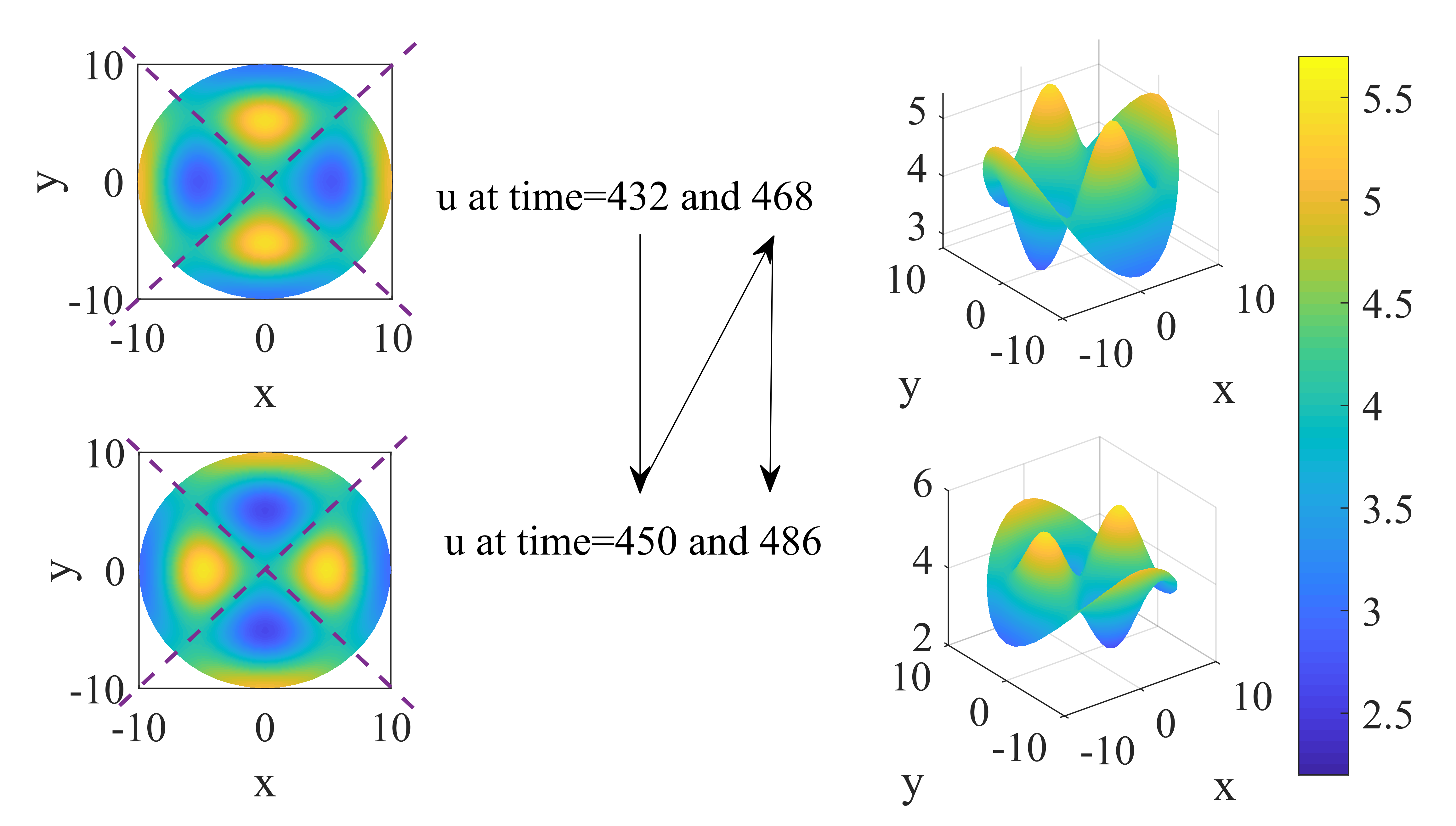}
(b)\includegraphics[width=0.47\textwidth]{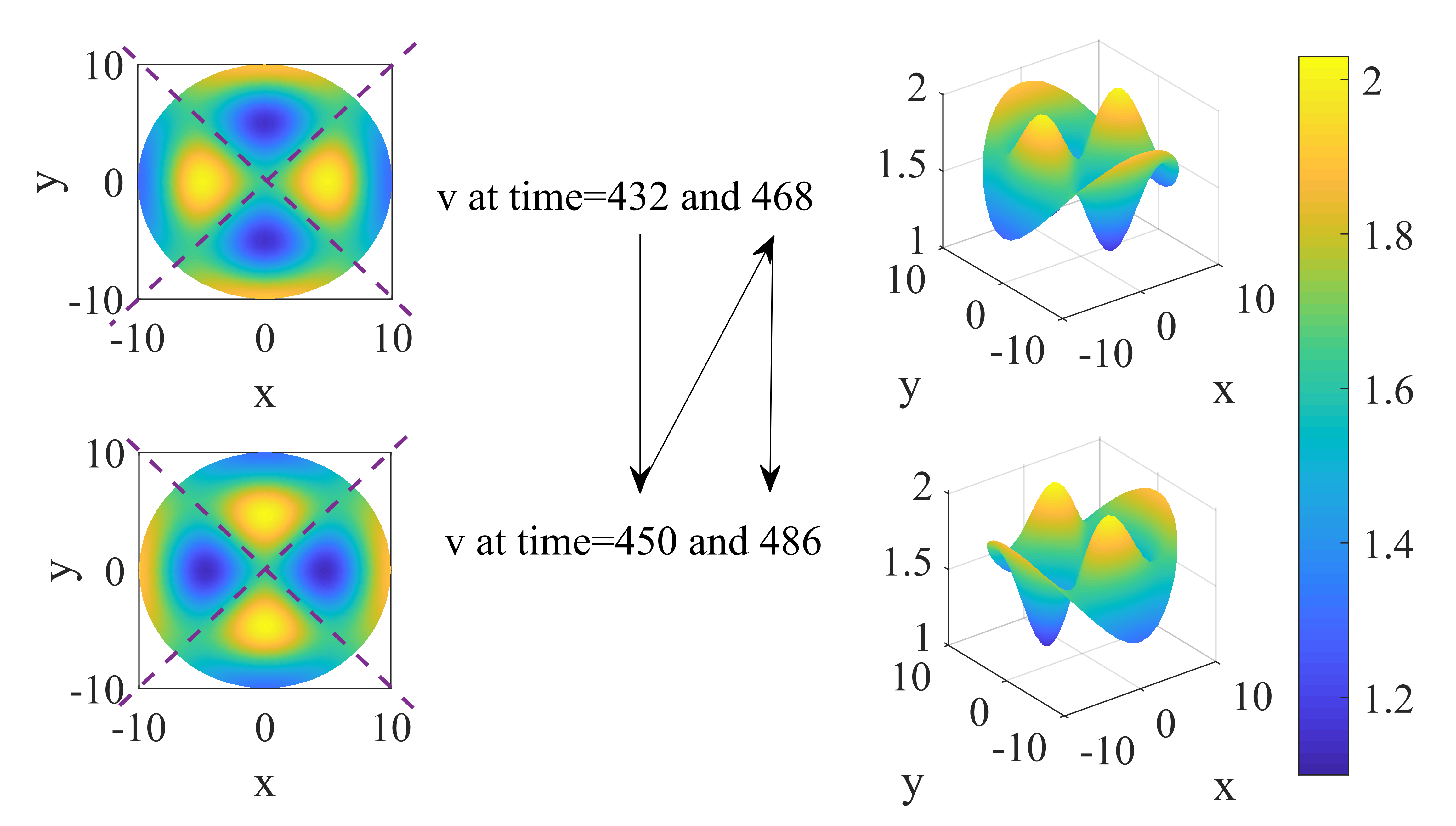}
\caption{System (\ref{system R-M}) generates standing wave that has two fixed axes with $(\chi,\tau)=(0.46,9.6)$.
Initial values are $u(t,r,\theta)=u^*(1+0.1\cdot\cos{t}\cdot\cos{\frac{2\pi r}{R}}\cdot\cos{2\theta}),~v(t,r,\theta)=v^*(1+0.1\cdot\cos{t}\cdot\cos{\frac{2\pi r}{R}}\cdot\cos{2\theta}),~t\in[-\tau,0)$. (a): u; (b): v.}
\label{S-1}
\end{figure}

\begin{figure}[t]
\centering
(a)\includegraphics[width=0.459\textwidth]{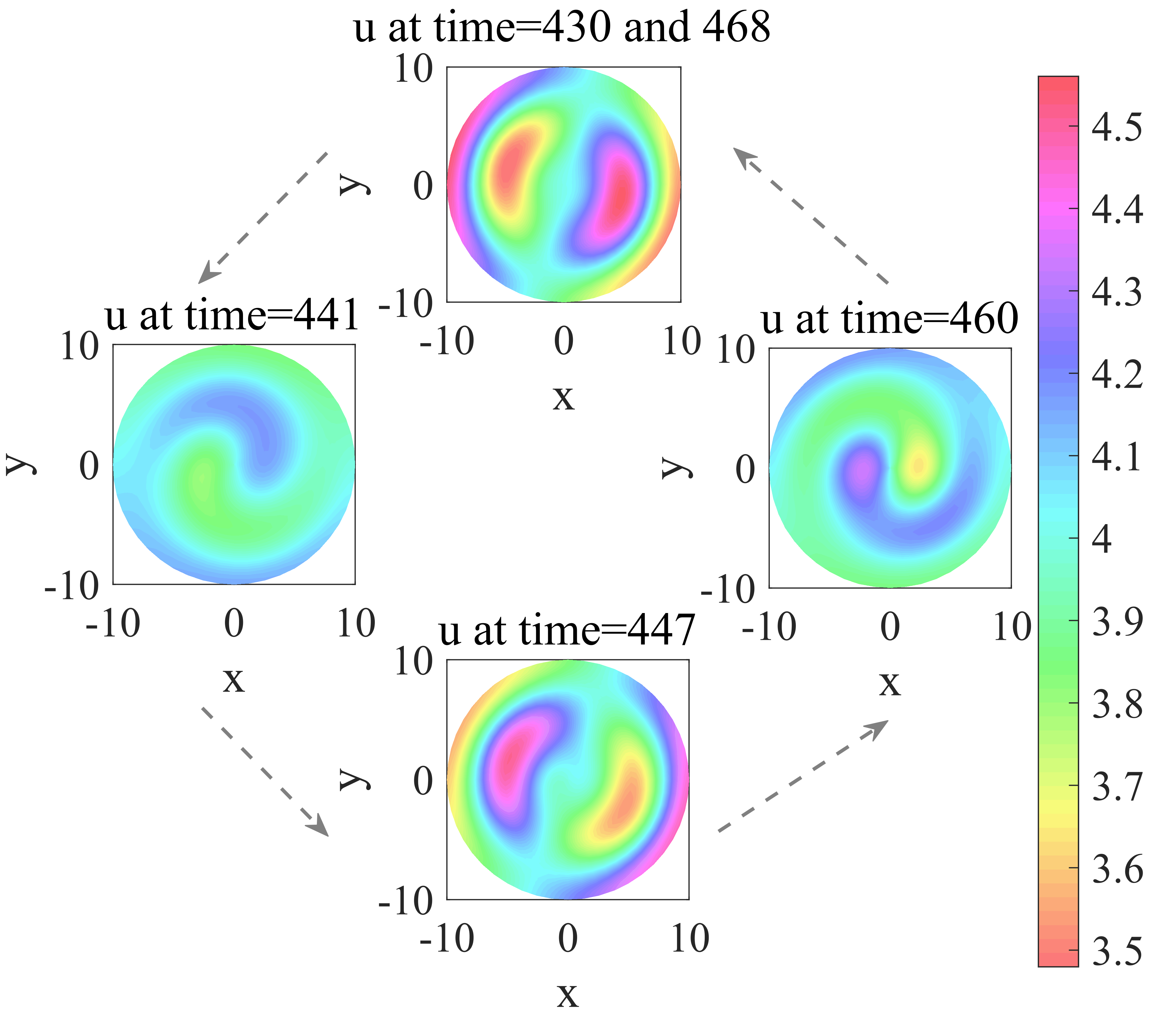}
(b)\includegraphics[width=0.459\textwidth]{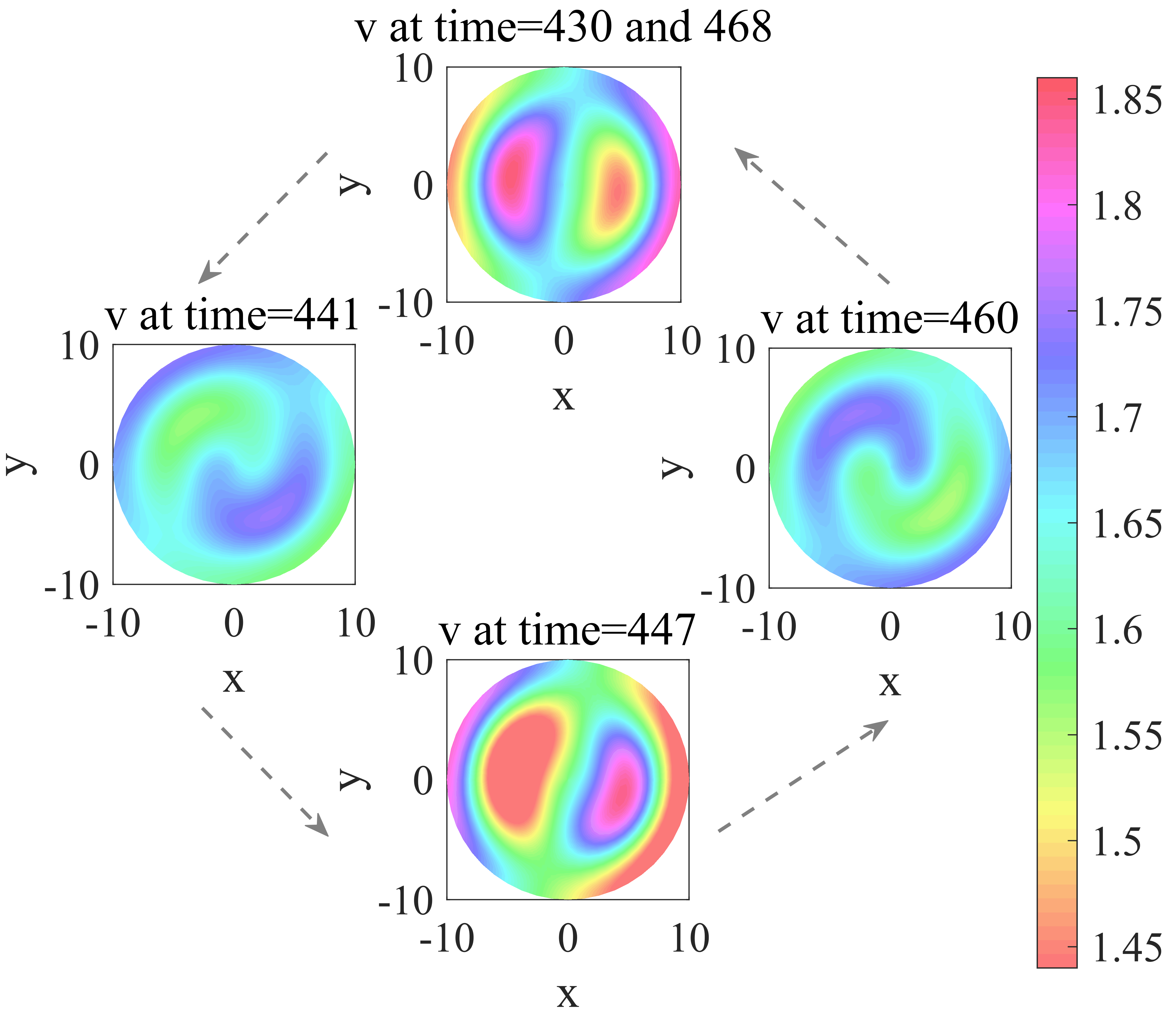}
\caption{
System (\ref{system R-M}) generates counterclockwise rotating wave with $(\chi,\tau)=(0.38,9.88)$. The patterns at four different time within a period are selected to show the periodic changes in population distribution.
Initial values are $u(t,r,\theta)=u^*(1+0.1\cdot\cos{t}\cdot\cos{\frac{2\pi r}{R}}\cdot\sin{\theta})$,
$v(t,r,\theta)=v^*(1+0.1\cdot\cos{t}\cdot\cos{\frac{2\pi r}{R}}\cdot\cos{\theta}),~t\in[-\tau,0)$. (a): u; (b): v.}
\label{R-2}
\end{figure}

\begin{figure}[t]
\centering
(a)\includegraphics[width=0.459\textwidth]{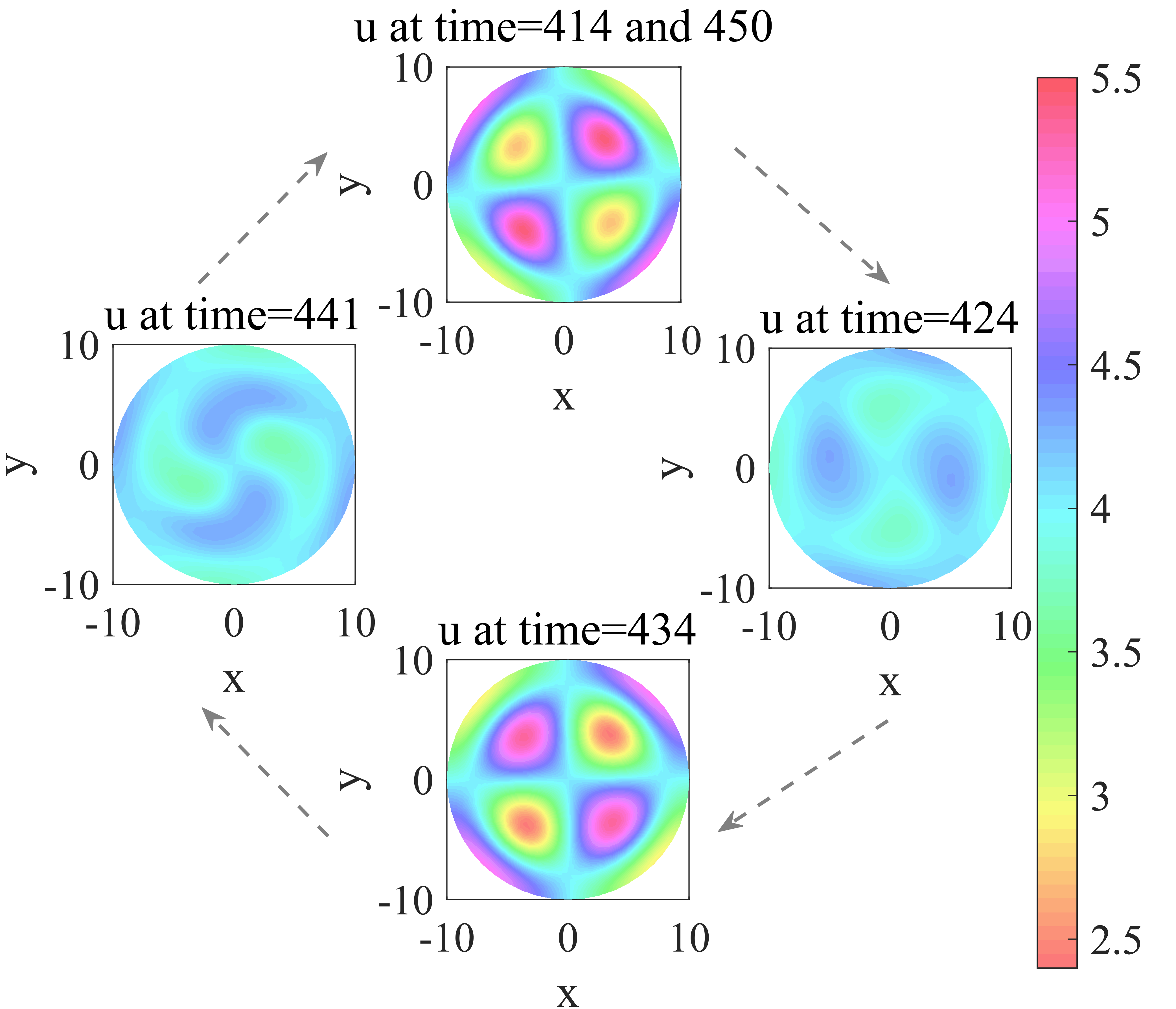}
(b)\includegraphics[width=0.459\textwidth]{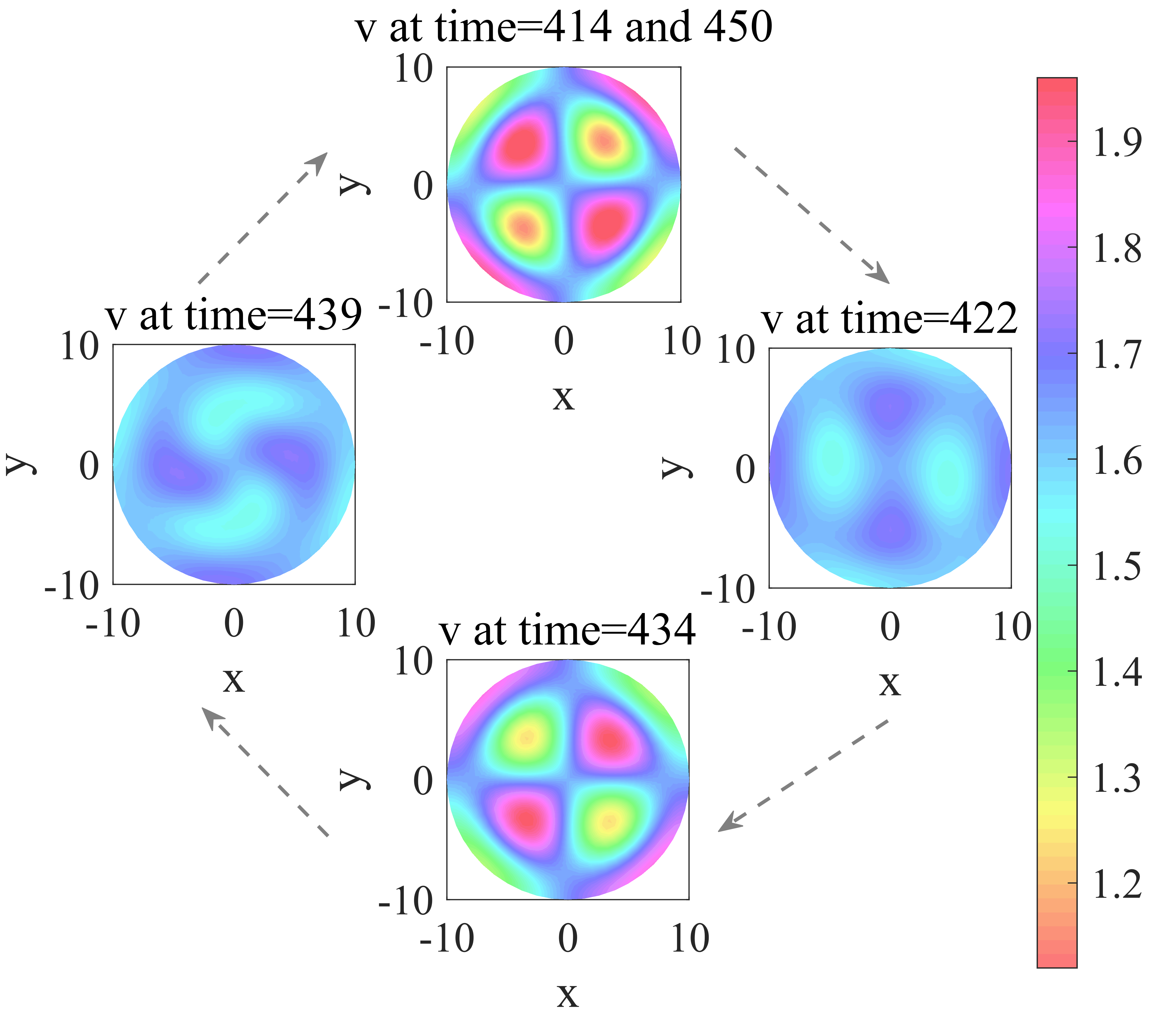}
\caption{
System (\ref{system R-M}) generates clockwise rotating wave with $(\chi,\tau)=(0.46,9.6)$. The patterns at four different time within a period are selected to show the periodic changes in population distribution.
Initial values are $u(t,r,\theta)=u^*(1+0.1\cdot\cos{t}\cdot\cos{\frac{2\pi r}{R}}\cdot\cos{2\theta}),~v(t,r,\theta)=v^*(1+0.1\cdot\cos{t}\cdot\cos{\frac{2\pi r}{R}}\cdot\sin{2\theta}),~t\in[-\tau,0)$. (a): u; (b): v.}
\label{R-1}
\end{figure}

\begin{figure}[t]
\centering
(a)\includegraphics[width=1\textwidth]{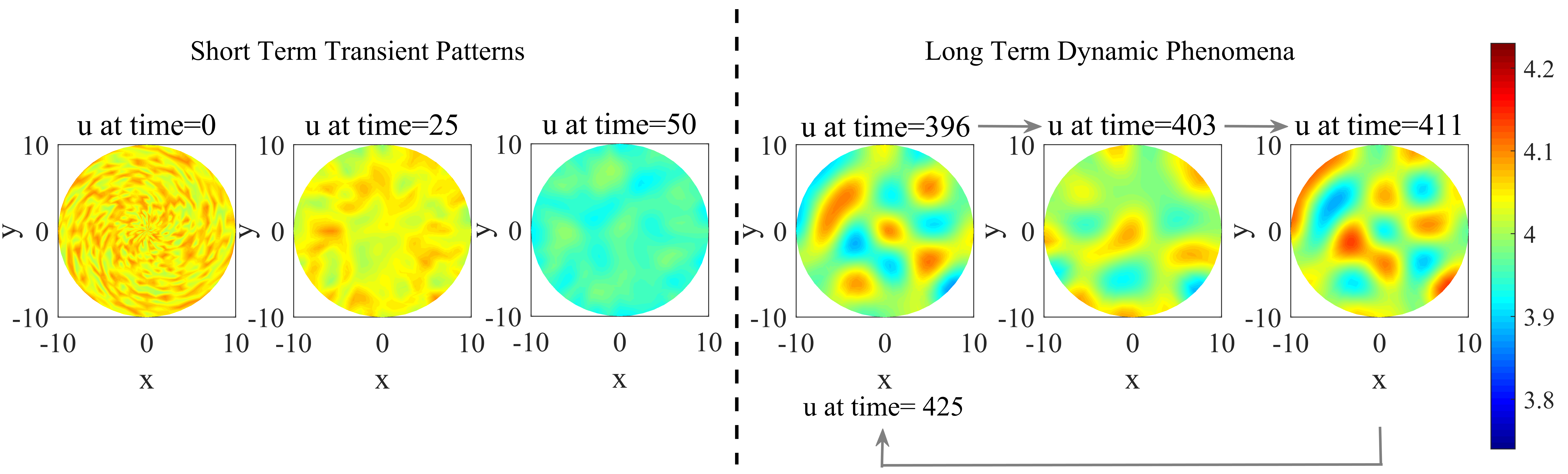}
(b)\includegraphics[width=1\textwidth]{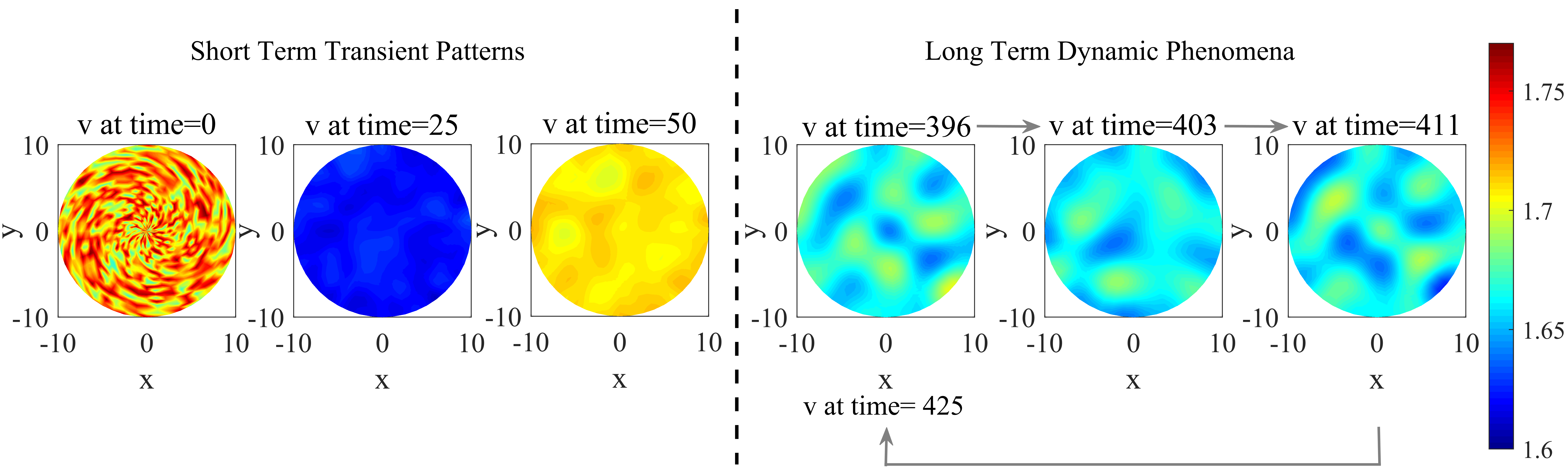}
\caption{System (\ref{system R-M}) generates rich transient patterns within a short term and periodic patterns with complex spatial structures after long-term evolution under random initial value and $(\chi,\tau)=(0.46,9.6)$. The patterns at three
different time within a period are selected to show the periodic changes in population distribution.
(a): u; (b): v.}
\label{rand}
\end{figure}

Fixing $d_1=0.1,~d_2=0.2,~\alpha=1,~K=6,~d=0.1,~R=10$, at the unique positive constant steady solution $(u^*,v^*)=(4,1.67)$, partial bifurcation curves on the $\chi-\tau$ plane are shown in Figure \ref{fenzhitutaxis}.
Selecting $(\chi,\tau)=(0.38,9.88)$ (Case 1) or $(\chi,\tau)=(0.46,9.6)$ (Case 2), the standing wave (see Figure \ref{S-2} and Figure \ref{S-1}), counterclockwise (see Figure \ref{R-2}) and clockwise rotating wave (see Figure \ref{R-1}) can be found. For simplicity, patterns of the clockwise rotating wave under Case 1 and the counterclockwise rotating wave under Case 2 have not been listed. Besides, more complex periodic spatiotemporal patterns (see Figure \ref{rand}) with random initial value can be found.

Through numerical calculations, we obtain that under Case 1, $\omega^*=0.1567,~\tau_{{11}}^0=9.8270.~{\tau}^{\prime}(0) \approx -0.4085$, ${\rho}^{\prime}(0)\approx 0.1390$ for $\psi=\mathrm{e}^{\mathrm{i}\omega^*\vartheta} V_1 \hat{\phi}_{11}^c$, and ${\tau}^{\prime}(0) \approx -0.0883,~{\rho}^{\prime}(0)\approx 0.0402$ for $\psi=\mathrm{e}^{\mathrm{i}\omega^*\vartheta} V_1 \hat{\phi}_{11}^c+\mathrm{e}^{\mathrm{i}\omega^*\vartheta} V_1 \hat{\phi}_{11}^s$. Under Case 2, $\omega^*=0.1938,~\tau_{{21}}^0=9.5520.~{\tau}^{\prime}(0) \approx -0.0631$, ${\rho}^{\prime}(0)\approx 0.0467$ for $\psi=\mathrm{e}^{\mathrm{i}\omega^*\vartheta} V_1 \hat{\phi}_{21}^s$, and ${\tau}^{\prime}(0) \approx -0.4942,~{\rho}^{\prime}(0)\approx 0.1697$ for $\psi=\mathrm{e}^{\mathrm{i}\omega^*\vartheta} V_1 \hat{\phi}_{21}^c+\mathrm{e}^{\mathrm{i}\omega^*\vartheta} V_1 \hat{\phi}_{21}^s$.
By Theorem \ref{xingzhiLS}, four bifurcation periodic solutions under two cases are both supercritical the period of them are both greater than $\frac{2\pi}{\omega_*}$.

\begin{remark}
On a one-dimensional interval, the solution often quickly converges to a $\cos{kx}$-like spatially inhomogeneous periodic solution \citep{Dong2023J,Shi2022J}. However, the two-dimensional domain yields a more significant process of aggregation and dispersion, ultimately presenting periodic patterns with complex spatial structures. The standing wave, rotating wave, and more complex spatiotemporal patterns under different initial conditions on a disk are new interesting phenomena.
\end{remark}

\section{Conclusion}

In this paper, we mainly investigated the equivariant Hopf bifurcation bifurcating from the positive equilibrium
of \eqref{system R-M}.  Methods of Lyapunov-Schmidt reduction and isotropic subgroups were combined to explore the interaction between symmetry, time delay, and taxis on a disk. We have done numerical simulations and spatially inhomogeneous periodic solutions were found, including standing waves, rotating waves, and more complex patterns.

 The time delay and taxis can induce the generation of spatially inhomogeneous periodic solutions, while specific forms of standing and rotating waves will be generated on the disk. It is worth mentioning that the algorithm in this paper can also be extended to other fields, such as chemistry, mechanics, nonlinear optics, etc. This extension allows for the exploration of the properties of standing and rotating waves, offering effective control in diverse applications.
 
 \bibliography{taxis}

\end{document}